\documentclass[12pt]{article}
\usepackage{amsmath, amsfonts, amssymb}
\newcommand{\R}{{\mathbf R}}
\newcommand{\ve}{\varepsilon}
\newcommand{\pa}{\partial}

\newtheorem{theorem}{Theorem}[section]
\newtheorem{proposition}[theorem]{Proposition}
\newtheorem{corollary}[theorem]{Corollary}
\newtheorem{lemma}[theorem]{Lemma}

\begin{document}
\title{Generalized wave operators for a system of nonlinear wave equations
in three space dimensions}

\author{Hideo Kubo \ (Corresponding author)
\\
Division of Mathematics,
Graduate School of Information Sciences,
\\
Tohoku University,
Sendai 980-8579, Japan
\\
Phone and Fax:\ +81-22-795-4628
\\
E-mail: kubo@math.is.tohoku.ac.jp
\\
\\
K\^oji Kubota
\\
A professor emeritus from Hokkaido
University, 
\\
Sapporo 060-0808, Japan}



\date{}

\maketitle

\begin{abstract}
This paper is concerned with the final value problem for a system of nonlinear wave equations.
The main issue is to solve the problem for the case where the nonlinearity is of a long range type.
By assuming that the solution is spherically symmetric, 
we shall show global solvability of the final value problem around a suitable final state, and hence
the generalized wave operator and long range scattering operator can be constructed.
\end{abstract}

\noindent
\small{Keywords\,:\,Final value problem, system of nonlinear wave equations, generalized wave operator,
long range scattering operator.}

\section{Introduction}
This paper is a continuation of a previous one \cite{16}, in which we have
studied the global existence and asymptotic behavior of solutions to the system
of semilinear wave equations\,:
\begin{equation}\label{1.1}
\left\{\begin{array}{ll}
\partial_t^2 u-\Delta u=|\partial_t v|^p \quad 
&\mbox{in} \quad  {\mathbb R}^3 \times {\mathbb R},\\
\partial_t^2 v -\Delta v =|\partial_t u|^q \quad 
&\mbox{in}\quad  {\mathbb R}^3 \times {\mathbb R},
\end{array}\right.
\end{equation}
where $1<p\le q$.
Assuming 
\begin{equation}\label{1.4}
 q(p-1)>2,
\end{equation}
we have shown that the initial value problem for (\ref{1.1}) has 
a global solution if the initial 
data are radilally symmetric and sufficiently small
(notice that if (\ref{1.4}) fails true, then the classical solution to 
(\ref{1.1}) blows up in finite time however small the initial data are,
in general, due to Deng \cite{ken2}). 
In the present paper we consider a problem whether the scattering operator 
for (\ref{1.1}) can be defined or not.

When $p>2$, the global solution has the \lq\lq free profile".
Therefore, in this case, one can expect that the scattering operator
is defined by $S=(W_+)^{-1} W_-$, following Lax and Phillips \cite{LaPh}.
Here the wave operators $W_+$ and $W_-$ are obtained by solving 
the final value problem for (\ref{1.1}). 
In fact, for a given final state $(u^+, v^+)$ which is a solution to the system 
of homogeneous wave equations\,:
\begin{equation}\label{1.5} 
\partial_t^2 u^+- \Delta u^+=0, \quad 
\partial_t^2 v^+- \Delta v^+=0 \quad 
\mbox{in} \quad  {\mathbb R}^3 \times {\mathbb R},
\end{equation}
if we found a unique solution $(u,v)$ 
in ${\mathbb R}^3 \times [0,\infty)$ to (\ref{1.1}) satisfying 
\begin{equation}\label{1.6}
\|u(t)-u^+(t)\|_{E}
+\|v(t)-v^+(t)\|_{E} \to 0 \quad \mbox{as} \quad t \to \infty,
\end{equation}
then $W_+$ is defined by 
$$
(u^+,\partial_t u^+,v^+,\partial_t v^+)(x,0) \longmapsto 
(u,\partial_t u, v, \partial_t v)(x,0).
$$
Here $\|w(t)\|_{E}$ stands for the energy norm of $w(x,t)$, i.e., 
$$
\|w(t)\|_{E}^2=\frac12 \int_{{\mathbb R}^3} (|\partial_t w(x,t)|^2+|\partial_x w(x,t)|^2) dx.
$$ 
Analogously, $W_-$ is defined by replacing \lq\lq\,$t \to \infty$" in (\ref{1.6})
by \lq\lq\,$t \to -\infty$".

On the other hand, when $1<p\le 2$, 
the nonlinearity becomes of a long range type in the sense that the solution to 
the initial value problem for \eqref{1.1} exists globally in time but does not approach to any free solution. 
In fact, we have shown in \cite{16}
that the energy of the global solution is not generically bounded
for large $t>0$, so that it can not be asymptotic 
to any free solution of finte energy.
Nevertheless, we proved that the global solution $(u,v)$ to the problem 
has a \lq\lq generalized profile". 
More precisely, letting $(\bar{w},\bar{v})$ be a solution of
\begin{equation}\label{1.7}
\partial_t^2 \bar{w}- \Delta \bar{w}=F(x,t), \quad 
\partial_t^2 \bar{v}- \Delta \bar{v}=0 
\quad \mbox{in} \quad  {\mathbb R}^3 \times {\mathbb R}
\end{equation}
with a suitable $F(x,t)$, we have
\begin{equation}\label{1.8}
\|u(t)-\bar{w}(t)\|_{E}
+\|v(t)-\bar{v}(t)\|_{E} \to 0 \quad \mbox{as} \quad t \to \infty.
\end{equation}
For instance, $F(x,t)=|\partial_t \bar{v}(x,t)|^p$ when $0 \le 2-p <q(p-1)-2$.
These results imply that it is impossible to construct the wave operators
$W_+$ and $W_-$, in general, and suggest us to modify the definition of 
wave operators for (\ref{1.1}) when $1<p\le 2$.   

As in \cite{16}, we deal with only radially symmetric solution 
of (\ref{1.1}) in the present paper. For this, we write
\begin{equation}\label{1.11}
u(x,t)=u_1(|x|,t), \quad v(x,t)=u_2(|x|,t).
\end{equation}
Then \eqref{1.1} becomes to 
\begin{equation}\label{1.12}
\left\{\begin{array}{ll}
\partial_t^2 u_1- \left(\partial_r^2+\frac{2}{r} \partial_r\right) 
u_1=|\partial_t u_2|^p \quad 
&\mbox{in} \quad  r>0, \ t \in {\mathbb R},\\
\partial_t^2 u_2- \left(\partial_r^2+\frac{2}{r} \partial_r\right) u_2
  =|\partial_t u_1|^q \quad 
&\mbox{in}\quad  r>0, \ t \in {\mathbb R}.
\end{array}\right.
\end{equation}
We wish to compare the nonlinear evolution under \eqref{1.12} with the 
linear evolution obeying
\begin{equation}\label{aux0.9}
\left\{\begin{array}{ll}
\partial_t^2 w- \left(\partial_r^2+\frac{2}{r} \partial_r\right) 
  w=F(r,t) \quad 
&\mbox{for} \quad  r>0, \ t \in {\mathbb R},\\
\partial_t^2 v- \left(\partial_r^2+\frac{2}{r} \partial_r\right) v
  =0 \quad 
&\mbox{for}\quad  r>0, \ t \in {\mathbb R}
\end{array}\right.
\end{equation}
with a suitably chosen $F(r,t)$, as $t \to \pm \infty$.
Actually, we are able to realize this by assuming
the following slightly stronger condition than \eqref{1.4}\,:
\begin{equation}\label{1.25} 
(p-1)^2 (q-1) >1
\end{equation}
(for the detail, see Theorems \ref{Th.2-} and \ref{Th.2} below).
Since we are able to solve the initial value problem in the same function
space (see Theorem \ref{Th.4} below), 
these results lead us to a construction of a long range scattering operator 
for (\ref{1.1}).

We remark that this kind of modification goes back to the seminal work of
Ozawa \cite{Oza91} for the nonlinear Schr\"odinger equation
(see also \cite{GinOza94, HayOza94, HayNau06, Wir07}, for instance).
To our knowledge, this paper provides the first result on the 
wave equation in this direction.
In order to treat the system \eqref{1.12}, we need to overcome 
a difficulty to handle the nonlinearities with small powers $p$
which can be close to $1$ under the assumption \eqref{1.25}.
For this reason, we shall construct a generalized final state,
which solves \eqref{aux0.9}, by using an iteration (see \eqref{1.25+}, \eqref{1.26+}).
Then we shall prove the solvability of the nonlinear system
\eqref{1.12} around the generalized final state
by introduing a suitable metric space
given by \eqref{kubota4.21}.

This paper is organized as follows.
In the next section we collect notation.
In the section 2 we present our main results.
The section 3 is a summary of \cite[Section 4]{16}.
We refine Theorem 6 and the part ({\rm ii}) of Theorem 7 in \cite{16}
so that one can take a parameter $\gamma$ to be positive.
The section 4 is devoted to prove the main theorems. 

\setcounter{equation}{0}
\section{Notation}



First we introduce a class of initial data\,:
\begin{eqnarray}\nonumber 
&& \quad Y_\nu(\varepsilon)=\{\vec{f}=(f,g) \in C^1({\mathbb R}) \times 
C({\mathbb R})~;~\vec{f}(-r)=\vec{f}(r) \ (r \in {\mathbb R}),
\\ \nonumber
&& \quad \quad \quad
r\vec{f}(r) \in C^2({\mathbb R}) \times C^1({\mathbb R}) \
 \mbox{and} \ \sup_{r>0}\,(1+r)^\nu 
|\!|\!| \vec{f}(r) |\!|\!| \le \varepsilon \},
\end{eqnarray}
where $\nu\in {\mathbb R}$, $\varepsilon>0$ and
\begin{eqnarray*}
|\!|\!| \vec{f}(r) |\!|\!|=|f(r)|+(1+r)(|f'(r)|+|g(r)|)+r(|f^{''}(r)|+|g'(r)|).
\end{eqnarray*}

Next we define several function spaces and norms.
Let $s=1$ or $s=2$.
First of all, we introduce a basic space of our argument\,:
\begin{eqnarray}\nonumber 
&& {X}^s=\{u(r,t) \in C^{s-1}({\mathbb R} \times [0,\infty))~;~
ru(r,t) \in C^s({\mathbb R} \times [0,\infty)),
\\ \nonumber
&& \hspace{15mm} \ \ u(-r,t)=u(r,t) \ \ \text{for}\ (r,t) \in {\mathbb R} \times [0,\infty) \}.
\end{eqnarray}
For $r>0$ and $t \ge 0$ we put
\begin{equation}\nonumber 
{[u(r,t)]_2}=|u(r,t)|+(1+r)\sum_{|\alpha|=1} |\partial^\alpha u(r,t)|+
r\sum_{|\alpha|=2} |\partial^\alpha u(r,t)|
\end{equation}
if $u \in X^2$, and 
\begin{equation}\nonumber 
{[u(r,t)]_1}=|u(r,t)|+r\sum_{|\alpha|=1} |\partial^\alpha u(r,t)|
\end{equation}
if $u \in X^1$, where  
$\partial=(\partial_r,\partial_t)$ and $\alpha$ is a multi-index.
For $\nu \in {\mathbb R}$, 
we define Banach spaces\,: 
\begin{eqnarray}\nonumber 
&& X^s(\nu)=\{u(r,t) \in X^s~;~ \|u\|_{X^s(\nu)} <\infty\},
\\ \nonumber
&& Z^s(\nu)=\{u(r,t) \in X^s~;~ \|u\|_{Z^s(\nu)} <\infty\},
\end{eqnarray}
where we have set
\begin{align}\label{2.8}
& \|u\|_{X^s(\nu)}=\sup_{r>0, \, t \ge 0}
[u(r,t)]_s\,(1+|r-t|)^{\nu}, 
\\ \label{2.9}
& \|u\|_{Z^s(\nu)}=\sup_{r>0, \, t \ge 0}
[u(r,t)]_s\,(1+r+t)^{\nu-1} (1+|r-t|).
\end{align}
Notice that $X^s(\nu) \subset Z^s(\nu)$ if $\nu \le 1$,
while $Z^s(\nu) \subset X^s(\nu)$ if $\nu \ge 1$.


For notational symplicity, we shall denote $\| w(|\cdot|,t)\|_{E}$ by 
$\| w(t)\|_{E}$ for a function $w(r,t)$.


\setcounter{equation}{0}
\section{Main Results}

\subsection{Existence of wave operators}

When $p>2$, the evolution obeying \eqref{1.12} is well characterized
by the homogeneous wave equation.
For this, we first recall the known fact about the initial value problem for the homogeneous 
wave equation (see e.g. \cite{16})\,:
\begin{eqnarray}\label{3.1}
&& u_{tt}-\left(u_{rr}+\frac{2}{r} u_r\right)=0
\quad\mbox{in}\quad (0,\infty)\times(0,\infty),
\\
\label{3.2}
&& 
(u, \partial_t u)(r,0)=\vec{f}(r) \ 
\quad\mbox{for}\quad r>0.
\end{eqnarray}
The solution of this problem is expressed by
\begin{equation}\label{2.1}
K[\vec{f}](r,t)=\frac1{2r} \left\{
\int_{r-t}^{r+t} \lambda g(\lambda) d\lambda
+\frac{\partial}{\partial t} \int_{r-t}^{r+t}
 \lambda f(\lambda) d\lambda \right\}.
\end{equation}
Moreover we have

\begin{proposition} \label{Th.3.1}
Let $\varepsilon>0$, $\nu>0$.
If $\vec{f} \in Y_\nu(\varepsilon)$, then $K[\vec{f}] \in X^2(\nu)$ and
\begin{equation}\label{3.3}
\|K[\vec{f}]\|_{X^2(\nu)}\le C \varepsilon
\end{equation}
holds, where $C$ is a constant depending only on $\nu$.
\end{proposition}

We set
\begin{eqnarray}\label{2.12}
\kappa_1=p-1, \quad 
\kappa_2=q-1 
\end{eqnarray}
for $p>2$. Then we see $\kappa_1$, $\kappa_2>1$.
Our main result in this subsection is as follows.

\begin{theorem}\label{Th.1}\,$(${\text Existence of a wave operator}$)$\,
Let $1<p \le q$.
Suppose that $p>2$. 
Then there is a positive number $\varepsilon_0$ $($depending only on
$p$ and $q)$ such that 
one can define the wave operator $W_+=(W_+^{(1)},W_+^{(2)})$ from
$Y_{\kappa_1}(\varepsilon_0) \times Y_{\kappa_2}(\varepsilon_0)$ to
$Y_{\kappa_1}(2\varepsilon_0) \times Y_{\kappa_2}(2\varepsilon_0)$ by
\begin{equation}\label{1.19}
W_+^{(j)}[\vec{f}_1,\vec{f}_2](r)=(u_j, \partial_t u_j)(r,0)
\quad (j=1,2),
\end{equation}
where $(u_1, u_2) \in X^2(\kappa_1) \times X^2(\kappa_2)$ is a unique solution
of \eqref{1.12} satisfying
\begin{equation}\label{1.20}
\|u_1(t)-K[\vec{f_1}](t)\|_{E}+\|u_2(t)-K[\vec{f_2}](t)\|_{E} \to 0 
  \quad \mbox{as} \quad t \to \infty
\end{equation} 
for each $(\vec{f}_1, \vec{f}_2) \in Y_{\kappa_1}(\varepsilon_0) \times Y_{\kappa_2}(\varepsilon_0)$.
Moreover, we have for $r>0$
\begin{eqnarray}\label{1.22}
&& |\!|\!| W_+^{(1)}[\vec{f}_1,\vec{f}_2](r)-\vec{f_1}(r) |\!|\!| (1+r)^{\kappa_1} \le  
   C \varepsilon^p,
\\ \label{1.23}
&& |\!|\!| W_+^{(2)}[\vec{f}_1,\vec{f}_2](r)-\vec{f_2}(r) |\!|\!| (1+r)^{\kappa_2} \le  
   C \varepsilon^q,
\end{eqnarray} 
provided $\vec{f}_j \in Y_{\kappa_j}(\varepsilon)$~$(j=1,2)$ and $0<\ve\le \ve_0$,
where $C$ is a constant depending only on $p$ and $q$.
\end{theorem}

Our next step is to construct the inverse of $W_+$,
based on the existence result given in Theorem 1 of \cite{16}
about the initial value problem for \eqref{1.12} with
\begin{equation}\label{1.33}
(u_1, \partial_t u_1)(r,0)=\vec{\varphi}_1(r), \ 
(u_2, \partial_t u_2)(r,0)=\vec{\varphi}_2(r)
 \quad 
\mbox{for} \quad  r>0.
\end{equation}
Let $(u_1,u_2) \in X^2(\kappa_1) \times {X^2(\kappa_2)}$ be the unique solution of 
the problem satisfying
\begin{eqnarray}\label{1.32}
\| u_1\|_{X^2(\kappa_1)}+ 
\| u_2\|_{X^2(\kappa_2)} \le 2C_0\varepsilon
\end{eqnarray} 
with $C_0$ the canstant in (\ref{3.3}).
Note that $(u_1,u_2)$ satisfies the following system
of integral equations\,:
\begin{equation}\label{2.31}
 u_1=K[\vec{\varphi}_1]+L(|\partial_t u_{2}|^p),
\quad u_2=K[\vec{\varphi}_2]+L(|\partial_t u_{1}|^q).
\end{equation}
(See e.g. \cite{16}).
Using the solution $(u_1,u_2)$, we define
\begin{align}\label{ko19}
w=u_1-R(|\partial_t u_2|^p), \quad
v=u_2-R(|\partial_t u_{1}|^q),
\end{align}
where $L$ and $R$ are the integral operators associated with
the inhomogeneous wave equation whose definition will be given
in \eqref{2.2} and \eqref{2.3} below, respectively.
If we set
\begin{align}\label{ko20}
 \vec{f}_1(r)=(w(r,0), \pa_t w(r,0)),
\quad
 \vec{f}_2(r)=(v(r,0), \pa_t v(r,0)) 
\end{align}
for $r>0$, then we see that 
\begin{align}\label{ko21}
w=K[\vec{f}_1], \quad v=K[\vec{f}_2].
\end{align}

Now we state the result for the inverse of ${W}_+$.

\begin{theorem}\label{Th.3}\,$(${\text Existence of the inverse of a wave operator}$)$\,
Let the assumptions of Theorem \ref{Th.1} hold.
Then there exists a positive number $\varepsilon_0$ $($depending only on
$p$ and $q)$ such that for any $\ve \in (0,\ve_0]$,
one can define $(W_+)^{-1}$ by
$(\vec{\varphi}_1, \vec{\varphi}_2)
\in Y_{\kappa_1}(\varepsilon) \times Y_{\kappa_2}(\varepsilon) \longmapsto
(\vec{f}_1, \vec{f}_2) \in
Y_{\kappa_1}(2\varepsilon) \times Y_{\kappa_2}(2\varepsilon)$ 
so that \eqref{1.20} is valid.
Here $(u_1,u_2)$ is the solution of \eqref{2.31}
satisfying \eqref{1.32}, and
$(\vec{f}_1,\vec{f}_2)$ is defined by \eqref{ko20}
for $(\vec{\varphi}_1, \vec{\varphi}_2)
\in Y_{\kappa_1}(\varepsilon) \times Y_{\kappa_2}(\varepsilon)$. 

Moreover, we have for $r>0$
\begin{eqnarray}\label{1.34}
&& |\!|\!| \vec{f_1}(r)- \vec{\varphi}_1(r) |\!|\!| (1+r)^{\kappa_1} \le  
   C \varepsilon^p,
\\ \label{1.35}
&& |\!|\!| \vec{f_2}(r)- \vec{\varphi}_2(r) |\!|\!| (1+r)^{\kappa_2} \le  
   C \varepsilon^q,
\end{eqnarray} 
provided $\vec{\varphi}_j \in Y_{\kappa_j}(\varepsilon)$~$(j=1,2)$ and $0<\ve\le \ve_0$,
where $C$ is a constant depending only on $p$ and $q$.
\end{theorem}

\noindent
{\bf Remark.}\ 
Now we are in a position to conclude the existence of a scattering
operator for (\ref{1.1}).
As we have constructed $W_+$ in Theorem
\ref{Th.1}, we obtain $W_-$ as well. 
Taking the range of $W_-$ to be included
by that of $W_+$, we are able to define
the scattering operator by $S=(W_+)^{-1} W_-$.

\subsection{Existence of generalized wave operators}
In this subsection we consider the case where $1<p \le 2$.
We set
\begin{eqnarray}\label{2.13}
\kappa_1=p-1,
\quad 
\kappa_2=q(p-1)-1 
\end{eqnarray}
for $1<p<2$.
While, when $p=2$, we take $\kappa_1$ and $\kappa_2$ in such a way that
\begin{equation}\label{aux6}
0<\kappa_1<1<\kappa_2<q-1, \quad q\kappa_1=\kappa_2+1.
\end{equation}
For instance, $\kappa_1=(q+2)/(2q)$, $\kappa_2=q/2$ satisfy the above conditions.
Note that $0<\kappa_1<1$ and $\kappa_2>1$ in both cases, by the assumption \eqref{1.4}.


First of all, we present a result for a special case of
Theorem \ref{Th.2} below, because it would make easy to 
recognize the statement for the general case.
Namely, we assume that $1<p<2$ and the following stronger condition on $p, q$ than
\eqref{1.25}\,:
\begin{equation}\label{4}  
\kappa_1 \kappa_2=(p-1)(q(p-1)-1) > 1.
\end{equation}
In order to have an analogue to Theorem \ref{Th.1}, we replace the final 
state $(K[\vec{f_1}], K[\vec{f_2}])$ by $(w_1, v_0) \in Z^2(\kappa_1) \times X^2(\kappa_2)$ 
which is the solution of the initial value problem for 
\begin{equation}\label{1.27}
\left\{\begin{array}{ll}
\partial_t^2 w_1- \left(\partial_r^2+\frac{2}{r} \partial_r\right) 
w_1=|\partial_t v_0|^p
 \quad 
&\mbox{for} \quad  r>0, \ t \in {\mathbb R},\\
\partial_t^2 v_0- \left(\partial_r^2+\frac{2}{r} \partial_r\right) v_0
  =0 \quad 
&\mbox{for}\quad  r>0, \ t \in {\mathbb R}
\end{array}\right.
\end{equation}
with 
\begin{equation}\label{9}  
(w_1, \partial_t w_1)(r,0)=\vec{f}_1(r), \quad
(v_0, \partial_t v_0)(r,0)=\vec{f}_2(r) 
 \quad 
\mbox{for} \quad  r>0.
\end{equation}
Actually, we have the following.

\begin{theorem}\label{Th.2-}\,$(${\text Existence of a generalized wave 
operator\,;\,a special case}$)$\,
Let $1<p \le q$.
Suppose that $1<p<2$ and \eqref{4} holds.
Then there exists a positive number $\varepsilon_0$ $($depending only on
$p$ and $q)$ such that 
one can define a generalized wave operator $\widetilde{W}_+=(\widetilde{W}_+^{(1)},
\widetilde{W}_+^{(2)})$ from
$Y_{\kappa_1}(\varepsilon_0) \times Y_{\kappa_2}(\varepsilon_0)$ to
$Y_{\kappa_1}(2\varepsilon_0) \times Y_{\kappa_2}(2\varepsilon_0)$ by 
\begin{equation}\label{*}  
\widetilde{W}_+^{(j)}[\vec{f}_1,\vec{f}_2](r)=(u_j, \partial_t u_j)(r,0)
\quad (j=1,2),
\end{equation}
where $(u_1, u_2) \in Z^2(\kappa_1) \times X^2(\kappa_2)$ is a unique solution
of \eqref{1.12} satisfying 
\begin{equation}\label{8}  
 \|u_1(t)-w_1(t)\|_{E} +\|u_2(t)-v_0(t)\|_{E}
\to 0 \quad \mbox{as} \quad t \to \infty
\end{equation} 
for each $(\vec{f}_1, \vec{f}_2) \in Y_{\kappa_1}(\varepsilon_0) \times Y_{\kappa_2}(\varepsilon_0)$.
Moreover, we have for $r>0$
\begin{equation}\label{10}  
 |\!|\!| \widetilde{W}_+^{(1)}[\vec{f}_1,\vec{f}_2](r)-\vec{f_1}(r) |\!|\!| 
(1+r)^{\kappa_1} \le  C \varepsilon^{1+(p-1)q} (1+r)^{-\kappa_1(\kappa_2-1)},
\end{equation}
and
\begin{equation}\label{11}  
 |\!|\!| \widetilde{W}_+^{(2)}[\vec{f}_1,\vec{f}_2](r)-\vec{f_2}(r) |\!|\!| (1+r)^{\kappa_2} \le  
   C \varepsilon^q,
\end{equation}
provided $\vec{f}_j \in Y_{\kappa_j}(\varepsilon)$~$(j=1,2)$ and $0<\ve\le \ve_0$,
where $C$ is a constant depending only on $p$ and $q$.
\end{theorem}

\noindent
{\bf Remark.}\ 
When $p=2$, we have only to assume $q>2$, instead of \eqref{4}.
In fact, if we replace the right hand side of \eqref{10} by
$C\ve^{1+q}(1+r)^{-(\kappa_2-\kappa_1)}$, then the conclusions
of Theorem \ref{Th.2-} remain valid
(for the needed modification of the proof, see the remark 
after the proof of Theorem \ref{Th.2}).

\medskip

Next we relax the condition \eqref{4} to 
\begin{equation}\label{3}  
\kappa_1 \kappa_2 > 1+\kappa_1^2-\kappa_1,
\end{equation} 
which is equivalent to \eqref{1.25}, while we shall keep $1<p<2$.
In the previous case, it suffices to iterate just once for getting 
$w_1$ as a final state for $u_1$.
However, in order to treat the general case, we need to iterate
many times for finding out a suitable final state for $u_1$.

First we define a sequence $\{a_j\}_{j=0}^\infty$ by $a_0=1$ and
\begin{equation}\label{6} 
a_{j+1}=\kappa_1 (a_j-1)+\kappa_2 \quad \mbox{for}\ j \ge 0,
\end{equation} 
explicitly we have
$$
a_j=\frac{\kappa_2-\kappa_1}{1-\kappa_1}-
     \frac{(\kappa_2-1)(\kappa_1)^j}{1-\kappa_1}
\quad \mbox{for}\ j \ge 0.
$$
Observe that $\{a_j\}_{j=0}^\infty$ is strictly increasing,
$a_j<(\kappa_2-\kappa_1)/(1-\kappa_1)$ for $j\ge 1$, and
$\lim_{j\to\infty}a_j=(\kappa_2-\kappa_1)/(1-\kappa_1)$.
Since $p<2$ and \eqref{3} yield
$$
a_0=1<\frac{1}{\kappa_1}<\frac{\kappa_2-\kappa_1}{1-\kappa_1}
=\lim_{j\to\infty}a_j,
$$
there exists a nonnegative integer $\ell$ such that
\begin{equation}\label{7}  
a_{\ell+1}>\frac{1}{\kappa_1}, \quad a_\ell \le \frac{1}{\kappa_1}.
\end{equation} 

Next we introduce a sequence $\{(w_j,v_j)\}_{j=0}^{\ell+1}$ as follows\,:
For $(\vec{f}_1, \vec{f}_2) \in Y_{\kappa_1}(\varepsilon) 
\times Y_{\kappa_2}(\varepsilon)$ we set
\begin{eqnarray*}
&& w_{1}=w_0+L(|\partial_t v_{0}|^p), \quad w_0=K[\vec{f}_1], \\
&& v_{1}=v_0+R(|\partial_t w_1|^q),  \quad v_0=K[\vec{f}_2].
\end{eqnarray*}
Moreover, we define
\begin{eqnarray}\label{1.25+}
&& w_{j+1}=w_{j}+L(|\partial_t v_{j}|^p-|\partial_t v_{j-1}|^p)), 
\\ \label{1.26+}
&& v_{j+1}=v_{j}+R(|\partial_t w_{j+1}|^q-|\partial_t w_{j}|^q)
\end{eqnarray}
for $1 \le j \le \ell$.
Here $L$ and $R$ are the integral operators given
by \eqref{2.2} and \eqref{2.3}, respectively.
Then we see that for $1 \le j \le \ell+1$,
\begin{equation}\label{1.27+}
\left\{\begin{array}{ll}
\partial_t^2 w_j- \left(\partial_r^2+\frac{2}{r} \partial_r\right) 
w_j=|\partial_t v_{j-1}|^p
 \quad 
&\mbox{for} \quad  r>0, \ t \in {\mathbb R},\\
\partial_t^2 v_j- \left(\partial_r^2+\frac{2}{r} \partial_r\right) v_j
  =|\pa_t w_j|^q \quad 
&\mbox{for}\quad  r>0, \ t \in {\mathbb R}
\end{array}\right.
\end{equation}
and $(w_j, \partial_t w_j)(r,0)=\vec{f}_1(r)$ for $r>0$.

Now, the following theorem shows that $w_{\ell+1}$ is a final state for $u_1$.

\begin{theorem}\label{Th.2}\,$(${\text Existence of a generalized wave operator}$)$\,
Let $1<p \le q$.
Suppose that 
$1<p<2$ and \eqref{3}.
Assume $\kappa_1 a_\ell <1$ in addition to \eqref{7}.
Then there exists a positive number $\varepsilon_0$ $($depending only on
$p$ and $q)$ such that
one can define a generalized wave operator $\widetilde{W}_+=(\widetilde{W}_+^{(1)},
\widetilde{W}_+^{(2)})$ from
$Y_{\kappa_1}(\varepsilon_0) \times Y_{\kappa_2}(\varepsilon_0)$ to
$Y_{\kappa_1}(2\varepsilon_0) \times Y_{\kappa_2}(2\varepsilon_0)$ by 
\eqref{*},
where $(u_1, u_2) \in Z^2(\kappa_1) \times X^2(\kappa_2)$ is a unique solution
of \eqref{1.12} satisfying 
\begin{equation}\label{12}  
 \|u_1(t)-w_{\ell+1}(t)\|_{E} +\|u_2(t)-v_0(t)\|_{E}
\to 0 \quad \mbox{as} \quad t \to \infty
\end{equation} 
for each $(\vec{f}_1, \vec{f}_2) \in Y_{\kappa_1}(\varepsilon_0) \times Y_{\kappa_2}(\varepsilon_0)$.
Moreover, for $r>0$, we have \eqref{11} and
\begin{equation}\label{13}  
|\!|\!| \widetilde{W}_+^{(1)}[\vec{f}_1,\vec{f}_2](r)-\vec{f_1}(r) |\!|\!| (1+r)^{\kappa_1} 
\le  C \varepsilon^{B_\ell} (1+r)^{-\kappa_1(a_{\ell+1}-1)},
\end{equation} 
provided $\vec{f}_j \in Y_{\kappa_j}(\varepsilon)$~$(j=1,2)$ and $0<\ve\le \ve_0$,
where $C$ is a constant depending only on $p$ and $q$.
Here we put 
$$
B_\ell=1+(p-1)(q+\ell(p+q-2)).
$$
\end{theorem}

\noindent
{\bf Remark.}\ 
If $1<p<2$ and \eqref{4} holds, then \eqref{3} is valid and 
$\kappa_1a_\ell<1$ is satisfied for $\ell=0$. 
Therefore, Theorem \ref{Th.2-} follows from Theorem \ref{Th.2}.

On the one hand, suppose $\kappa_1a_\ell=1$ (notice that we have
$\ell \ge 1$ in this case).
Then we need to modify the statement of Theorem \ref{Th.2} a little.
Letting $\delta$ be a number satisfying
\begin{align}\label{aux7}
0<\delta<a_\ell-a_{\ell-1}, \quad \kappa_1^2\delta<\kappa_1a_{\ell+1}-1,
\end{align}
we define
\begin{align}\label{aux8}
a_\ell^\prime=a_{\ell}-\delta, \quad \text{and}\quad 
a_{\ell+1}^\prime=\kappa_1(a_{\ell}^\prime-1)+\kappa_2\,(=a_{\ell+1}-\kappa_1\delta).
\end{align}
Observing that
\begin{align}\label{aux9}
a_{\ell-1}<a_\ell^\prime<a_{\ell+1}^\prime, \quad 
\kappa_1 a_{\ell}^\prime<1, \quad\text{and}\quad \kappa_1 a_{\ell+1}^\prime>1,
\end{align}
we can show the statement of the theorem 
with $a_{\ell+1}$ in \eqref{13} replaced by $a_{\ell+1}^\prime$.

\vspace{2mm}

Our next step is to construct the inverse of $\widetilde{W}_+$,
based on the existence result in Theorem 1 of \cite{16}
for the initial value problem \eqref{1.12} and \eqref{1.33}.
Let $(u_1,u_2) \in Z^2(\kappa_1) \times {X^2(\kappa_2)}$ be the unique solution of 
\eqref{2.31} satisfying
\begin{eqnarray}\label{1.36}
\| u_1\|_{Z^2(\kappa_1)}+ 
\| u_2\|_{X^2(\kappa_2)} \le 2C_0\varepsilon
\end{eqnarray} 
with $C_0$ the canstant in (\ref{3.3}).
Using the solution, we set $w_0^*=K[\vec{\varphi}_1]$,
$v_0^*=u_2-R(|\partial_t u_{1}|^q)$.
Moreover, when $\ell \ge 1$, we define for $1 \le j \le \ell$
\begin{align}\label{2.32}
& w_j^*=w_0^*+L(|\partial_t v_{j-1}^*|^p),
\\ \label{2.33}
& v_j^*=v_0^*+R(|\partial_t w_{j}^*|^q).
\end{align}
We furter define 
\begin{align}\label{2.34}
 w^*=u_1-R(|\partial_t u_2|^p-|\partial_t v_{\ell}^*|^p),
\end{align}
which we wish to regard as a final state for $u_1$.
If we set
\begin{align}\label{2.35}
& \vec{f}_1(r)=(w^*(r,0), \pa_t w^*(r,0)) \quad \mbox{for}\quad  r>0,
\\ \label{2.36}
& \vec{f}_2(r)=(v_0^*(r,0), \pa_t v_0^*(r,0)) \quad \mbox{for}\quad  r>0,
\end{align}
then we see that $v_0^*$ and $w^*$ are represented as 
\begin{align}\label{2.38}
w^*=K[\vec{f}_1]+L(|\partial_t v_{\ell}^*|^p), \quad
 v_0^*=K[\vec{f}_2].
\end{align}

Now we state the result for the inverse of $\widetilde{W}_+$.

\begin{theorem}\label{Th.4}\,$(${\text Existence of the inverse of 
a generalized wave operator}$)$\,
Let the assumptions of Theorem \ref{Th.2} be fulfilled.
Then there exists a positive number $\varepsilon_0$ $($depending only on
$p$ and $q)$ such that, for any $\varepsilon \in (0,\varepsilon_0]$, 
one can define $(\widetilde{W}_+)^{-1}$ by
$(\vec{\varphi}_1, \vec{\varphi}_2)
\in Y_{\kappa_1}(\varepsilon) \times Y_{\kappa_2}(\varepsilon) \longmapsto
(\vec{f}_1, \vec{f}_2) \in
Y_{\kappa_1}(2\varepsilon) \times Y_{\kappa_2}(2\varepsilon)$ 
so that 
\begin{align}\label{2.40}
& \|u_1(t)-(K[\vec{f}_1]+L(|\partial_t v_{\ell}^*|^p))(t)\|_{E} \to 0 
\quad & \mbox{as} \quad t \to \infty
\\
\label{2.39}
& \|u_2(t)-K[\vec{f}_2](t)\|_{E} \to 0 
\quad & \mbox{as} \quad t \to \infty,
\end{align} 
hold.
Here $(u_1,u_2)$ is the solution of \eqref{2.31}
satisfying \eqref{1.36},
$(\vec{f}_1,\vec{f}_2)$ is defined by \eqref{2.35}, \eqref{2.36},
and $v_\ell^*$ is given by \eqref{2.33}
for $(\vec{\varphi}_1, \vec{\varphi}_2)
\in Y_{\kappa_1}(\varepsilon) \times Y_{\kappa_2}(\varepsilon)$. 

Moreover, for $r>0$, we have \eqref{1.35} and
\begin{equation}\label{1.39}
 |\!|\!| \vec{f_1}(r)- \vec{\varphi}_1(r) |\!|\!| (1+r)^{\kappa_1} 
 \le  C \varepsilon^{B_\ell} (1+r)^{-\kappa_1(a_{\ell+1}-1)},
\end{equation} 
provided $\vec{\varphi}_j \in Y_{\kappa_j}(\varepsilon)$~$(j=1,2)$ and $0<\ve\le \ve_0$,
where $C$ is a constant depending only on $p$ and $q$.
\end{theorem}

\noindent
{\bf Remark.}\ 
Now we are in a position to conclude the existence of a scattering
operator for (\ref{1.1}).
As we have constructed $\widetilde{W}_+$ in Theorem
\ref{Th.2}, we obtain $\widetilde{W}_-$ 
as well. 
Taking the range of $\widetilde{W}_-$ to be included
by that of $\widetilde{W}_+$, we are able to define 
the long range scattering operator by 
$\widetilde{S}=(\widetilde{W}_+)^{-1} \widetilde{W}_-$.

\setcounter{equation}{0}
\section{Inhomogeneous wave equations}
In this section we summarize the results of the section 4 in \cite{16}
for the case $a=c=1$.
The first one is concerned with the initial value problem for the inhomogeneous 
wave equation with the zero initial data:
\begin{eqnarray}\label{3.4}
&& u_{tt}-\left(u_{rr}+\frac{2}{r} u_r\right)=F(r,t)
\quad\mbox{in}\quad (0,\infty)\times(0,\infty),
\\
\label{3.5}
&& u(r,0)=(\partial_t u)(r,0)=0
\quad\mbox{for}\quad r>0.
\end{eqnarray}
The solution of this problem is given by\begin{equation}
L(F)(r,t)=\frac1{2r} \int_{0}^t ds \int_{r-(t-s)}^{r+(t-s)}
 \lambda F(\lambda,s) d\lambda.
\label{2.2}
\end{equation}
In order to study the qualitative property of $L(F)$, we set
\begin{eqnarray}
&& \quad M_0(F)=\sup_{r>0,\, t \ge 0} |F(r,t)| r^\alpha (1+r)^\beta 
(1+r+t)^\gamma (1+|r-t|)^\delta,
\label{2.18}
\\ \label{2.19}
&& \quad M_1(F)=M_0(F)
\\ \nonumber
&& \quad \quad +\sup_{r>0,\, t \ge 0} |\partial_r F(r,t)| r^{\alpha+1} (1+r)^{\beta-1} 
 (1+r+t)^\gamma (1+|r-t|)^{\delta}.
\end{eqnarray}
for $\alpha$, $\beta$, $\gamma$, and $\delta \in {\mathbb R}$.
Then we have

\begin{proposition} \label{Th.3.2} 
Let $F \in {X}^1$. Then we have $L(F) \in X^2$. 
Moreover, if $M_{s-1}(F)$ with $s=1$ or $s=2$ is finite for $\alpha<3-s$,
$\beta\in \R$, $\gamma\ge 0$, 
and $\delta>1$, then there exists a constant ${C}$ depending only on
$\alpha$, $\beta$, $\gamma$, and $\delta$ such that
\begin{align} \label{3.7}
& \| L(F)\|_{X^s(\nu)} \le {C} M_{s-1}(F)
\quad & \mbox{if}\quad \alpha+\beta+\gamma>2,
\\ \label{3.7+}
& \| L(F)\|_{Z^s(\alpha+\beta+\gamma-1)} \le {C} M_{s-1}(F)
\quad & \mbox{if}\quad 1<\alpha+\beta+\gamma<2,
\end{align}
where $\nu=\min(\alpha+\beta+\gamma-1,\delta)$.
\end{proposition}

\noindent{\it Proof.}\ 
Note that the statement follows from the case $\gamma=0$,
since $(1+\lambda+s)^{-\gamma} \le (1+\lambda)^{-\gamma}$ when $\gamma>0$.
Therefore, applying Theorem 6 in \cite{16} where the case $\gamma=0$ was shown,
we conclude the proof.
\hfill$\square$

Next we study the operator 
\begin{equation}
R(F)(r,t)=\frac1{2r} \int_t^\infty ds \int_{(s-t)-r}^{(s-t)+r}
 \lambda F(\lambda,s) d\lambda,
\label{2.3}
\end{equation}
related to the final value problem.
Indeed, if $F \in X^1$ and 
\begin{equation}
\sup_{r\ge 1,\, t \ge 0} |F(r,t)| (1+r)^\beta (1+r+t)^\gamma (1+|r-t|)^\delta < \infty
\label{2.20}
\end{equation}
for $\beta+\gamma>2$, $\delta \in \R$,
then we have $R(F) \in {X}^1$ and it satisfies the inhomogeneous wave equation\,{\rm :}
\begin{eqnarray}\label{3.19}
(\partial_t^2-\Delta)R(F)(|x|,t)=F(|x|,t)
\end{eqnarray}
in the distributional sense on ${\mathbf R}^3 \times (0,\infty)$.
The following result, which is a refinement of Theorem 7 in \cite{16}
in the sense that one can take a parameter $\gamma$ to be positive,
will play an essential role in this paper. 

\begin{proposition} \label{Th.3.3}
If $F \in X^1$ and $M_{s-1}(F)$ with $s=1$ or $s=2$ is finite for $\alpha<3-s$,
$\beta$, $\gamma \in \R$,  
and $\delta>1$ satisfying $\alpha+\beta+\gamma>2$,
then $R(F) \in X^s$ and there exists a constant ${C}$ depending only on
$\alpha$, $\beta$, $\gamma$, and $\delta$ such that
\begin{equation}\label{3.16}
\| R(F) \|_{Z^s(\mu+\gamma)} \le {C} M_{s-1}(F),
\end{equation}
where $\mu=\min(\alpha+\beta-1,\delta)$.
\end{proposition}

\noindent{\it Proof.}\ 
Since the statement for $\gamma\le 0$ was shown in Theorem 7 in \cite{16}, 
it suffices to prove it for $\gamma>0$.
Seeing the proof, we find that $R(F) \in X^s$ is valid also for $\gamma>0$.
Hence it remains to show \eqref{3.16}.

It follows from (\ref{2.3}) that 
\begin{equation}\nonumber
R(F)(r,t)=\frac1{2r} \int_t^\infty ds \int_{|(s-t)-r|}^{(s-t)+r}
 \lambda F(\lambda,s) d\lambda,
\end{equation}
since $\lambda F(\lambda,s)$ is odd in $\lambda$.
Observe that if $\lambda \ge {|(s-t)-r|}$ and $s \ge t$, then we have
$\lambda+s \ge r+t$,
so that $(1+\lambda+s)^{-\gamma} \le (1+r+t)^{-\gamma}$ when $\gamma>0$. 
Therefore, if $\alpha+\beta>2$, one can reduce the proof to the case 
$\gamma=0$ which was already shown in \cite{16}.

Suppose, on the contrary, that $\alpha+\beta \le 2$.
We take a positive number $\rho>0$ satisfying
\begin{equation}\label{3.21}
2-(\alpha+\beta)<\rho \le \delta+1-(\alpha+\beta), \quad
\rho<\gamma  
\end{equation}
and set $\beta^\prime=\beta+\rho$, $\gamma^\prime=\gamma-\rho$.
Then we have $\alpha+\beta^\prime>2$, $\gamma^\prime>0$ and $\delta>1$.
Applying the result in the preceding case with $\beta$ and $\gamma$
replaced by $\beta^\prime$ and $\gamma^\prime$ respectively,
we obtain the needed conclusion, because $\min(\alpha+\beta^\prime-1,\delta)
=\min(\alpha+\beta-1+\rho,\delta)=\alpha+\beta-1+\rho$
and $\mu=\min(\alpha+\beta-1,\delta)=\alpha+\beta-1$.
This completes the proof. 
\hfill$\square$

\setcounter{equation}{0}
\section{Proof of Main Results}

\subsection{Proof of Theorems \ref{Th.1} and \ref{Th.3}}

First we prove Theorem \ref{Th.1}. Suppose $p>2$.
Let $(\vec{f}_1, \vec{f}_2) \in
Y_{\kappa_1}(\varepsilon) \times Y_{\kappa_2}(\varepsilon)$
with $0<\ve \le 1$, and set $w_0=K[\vec{f}_1]$, $v_0=K[\vec{f}_2]$.
Then it follows from \eqref{3.3} that
\begin{align} \label{ko0}
\|w_0 \|_{X^2(\kappa_1)}+\|v_0 \|_{X^2(\kappa_2)} \le C\ve.
\end{align}
Recall that $\kappa_1=p-1>1$ and $\kappa_2=q-1>1$.

We shall solve the following system of integral equations\,:
\begin{align}\label{ko1} 
 u_1=w_{0}+R(|\pa_t u_2|^p), \quad
 u_2=v_{0}+R(|\pa_t u_1|^q),
\end{align}
where $R$ is defined by \eqref{2.3}.
To this end, we define $T(u_1,u_2)=(T^{(1)}(u_2),T^{(2)}(u_1))$ by
\begin{align}\label{ko2} 
T^{(1)}(u_2)=w_{0}+R(|\pa_t u_2|^p),
\quad T^{(2)}(u_1)=v_{0}+R(|\pa_t u_1|^q).
\end{align} 
For $\ve>0$ we introduce a metric space
\begin{align} \label{ko3} 
 D_\ve=\{(u_1,u_2) \in X^2 \times X^2\,;\,d((u_1,u_2),(w_{0},v_{0}))
   \le \ve \},
\end{align}
where we have set 
$$
d((u_1,u_2),(u_1^*,u_2^*))=\| u_1-u_1^*\|_{Z^2(\kappa_1)}
 +\| u_2-u_2^*\|_{Z^2(\kappa_2)}.
$$

First we prepare the following.

\begin{lemma} \label{Lemmako1}
Let $(u_1,u_2) \in D_\ve$.
Then we have
\begin{align}\label{ko6} 
& \|T^{(1)}(u_2)-w_0\|_{Z^2(\kappa_1)} \le 
 C \ve^{p}, 
\\ \label{ko7}
& \|T^{(2)}(u_1)-v_0\|_{Z^2(\kappa_2)} \le 
 C \ve^{q}.
\end{align}
Moreover we have
\begin{align}\label{ko10} 
& \|T^{(1)}(u_2)-T^{(1)}(u_2^*)\|_{Z^2(\kappa_1)} \le 
 C \ve^{p-1} \|u_2-u_2^*\|_{Z^2(\kappa_2)}, 
\\ \label{ko11}
& \|T^{(2)}(u_1)-T^{(2)}(u_1^*)\|_{Z^2(\kappa_2)} \le 
 C \ve^{q-1} \|u_1-u_1^*\|_{Z^2(\kappa_1)}
\end{align}
for $(u_1,u_2)$, $(u_1^*,u_2^*) \in D_\ve$.
\end{lemma}

\noindent{\it Proof.}\ 
First we observe that if $(u_1,u_2) \in D_\ve$,
then we have
\begin{align} \label{ko5}
\|u_1\|_{X^2(\kappa_1)}+\|u_2\|_{X^2(\kappa_2)} \le C\ve,
\end{align}
due to $\kappa_1$, $\kappa_2>1$ and \eqref{ko0}.

We start with the proof of \eqref{ko6}.
In view of \eqref{ko2}, it suffices to show
\begin{align} \label{ko9}
\|R(|\pa_t u_2|^p)\|_{Z^2(\kappa_1)} \le C\ve^p.
\end{align}
We see from \eqref{ko5} that
$M_1(|\pa_t u_2|^p) \le C\ve^p$ holds
for $\alpha=\gamma=0$, $\beta=p$ and $\delta=p\kappa_2$,
where $M_1(F)$ is defined by \eqref{2.18} and \eqref{2.19}.
Since $\alpha+\beta+\gamma-1=\kappa_1>1$,
by \eqref{3.16} with $s=2$ we get \eqref{ko9},
which implies \eqref{ko6}.
Analogously we obtain \eqref{ko7}, because $\kappa_2>1$.

Next we show \eqref{ko10}.
It follows from \eqref{ko2} that
\begin{align} \label{ko12}
T^{(1)}(u_2)-T^{(1)}(u_2^*)=R(|\pa_t u_2|^p-|\pa_t u_2^*|^p).
\end{align}
Since $p>2$, we see from \eqref{ko5} that
$$
M_1(|\pa_t u_2|^p-|\pa_t u_2^*|^p) \le C
 \ve^{p-1} \|u_2-u_2^*\|_{Z^2(\kappa_2)}
$$
for $\alpha=0$, $\beta=p$, $\gamma=\kappa_2-1$, 
and $\delta=1+(p-1)\kappa_2$.
Since $\alpha+\beta+\gamma-1=\kappa_1+\kappa_2-1>1$,
by \eqref{3.16} with $s=2$ we obtain
\begin{align*}
 \|R(|\pa_t u_2|^p-|\pa_t u_2^*|^p)\|_{Z^2(\kappa_1+\kappa_2-1)} \le 
 C \ve^{p-1} \|u_2-u_2^*\|_{Z^2(\kappa_2)}.
\end{align*}
In view of \eqref{ko12}, we get \eqref{ko10}, since $\kappa_2>1$.
Analogously we have \eqref{ko11}, because $q \ge p>2$.
This completes the proof.
\hfill$\square$

\noindent
{\it End of the proof of Theorem \ref{Th.1}.}\
We see from Lemma \ref{Lemmako1} that
there exists a positive number $\ve_0$ depending only on
$p$ and $q$ such that if $0<\ve \le \ve_0$,
then we have $T(u_1,u_2) \in D_\ve$ and 
$$
d(T(u_1,u_2),T(u_1^*,u_2^*)) \le 2^{-1} 
d((u_1,u_2),(u_1^*,u_2^*))
$$ 
for $(u_1,u_2)$, $(u_1^*,u_2^*) \in D_\ve$, 
namely, $T$ is a contaction on $D_\ve$.
Hence we find a unique solution  
$(u_1,u_2) \in D_\ve$ of \eqref{ko1}.
Here and in what follows, we suppose that $0<\ve \le \ve_0$
and $(u_1,u_2)$ is the solution.

Since $T^{(1)}(u_2)=u_1$ and $w_0=K[\vec{f}_1]$,
it follows from \eqref{ko6} that
\begin{align} \label{ko15}
[(u_1-K[\vec{f}_1])(r,t)]_2 \le C\ve^p (1+r+t)^{-(\kappa_1-1)}
 (1+|r-t|)^{-1}.
\end{align}
Therefore we have
\begin{align} \label{ko13}
& \|(u_1-K[\vec{f}_1])(t)\|_E
\\ \nonumber \le
 & C\ve^p (1+t)^{-(\kappa_1-1)}
  \left(\int_0^\infty (1+|r-t|)^{-2} dr \right)^{1/2}
\\ \nonumber
 \le & C\ve^p (1+t)^{-(\kappa_1-1)}
\end{align}
for $t \ge 0$. Analogously by \eqref{ko7} we get 
\begin{align} \label{ko14}
 \|(u_2-K[\vec{f}_2])(t)\|_E
   \le C\ve^q (1+t)^{-(\kappa_2-1)}
\end{align}
for $t \ge 0$. Hence we obtain \eqref{1.20}.

Moreover, we easily get \eqref{1.22} by taking $t=0$ in \eqref{ko15}.
Analogously, \eqref{1.23} follows from \eqref{ko7}.
Thus we prove Theorem \ref{Th.1}.
\hfill $\square$ 
 
\vspace{2mm}

Next we show Theorem \ref{Th.3}.
Let $(u_1,u_2) \in X^2(\kappa_1) \times {X^2(\kappa_2)}$ be the unique solution of 
\eqref{2.31} satisfying \eqref{1.32}.
Then we see from \eqref{ko19} and \eqref{ko21} that
\eqref{1.20}, \eqref{1.34} and \eqref{1.35} follows from
\eqref{ko9} and
\begin{align} \label{ko17}
\|R(|\pa_t u_1|^q)\|_{Z^2(\kappa_2)} \le C\ve^q.
\end{align}
By virtue of \eqref{1.32}, \eqref{ko9} can be shown as before.
In the same way we obtain \eqref{ko17}.
Thus we prove Theorem \ref{Th.3}.
\hfill $\square$

\subsection{Proof of Theorems \ref{Th.2} and \ref{Th.4}}

We start by showing the following basic estimates. 

\begin{lemma} \label{Lemma1}
Let $v \in X^2(\kappa_2)$.
Then $L(|\pa_t v|^p) \in Z^2(\kappa_1)$ and we have
\begin{align} \label{ta1}
\| L(|\pa_t v|^p) \|_{Z^2(\kappa_1)} \le
 C \|v\|_{X^2(\kappa_2)}^p.
\end{align}

While, let $w \in Z^2(\kappa_1)$.
Then $R(|\pa_t w|^q) \in Z^2(\kappa_2)$ and we have
\begin{align} \label{ta2}
\| R(|\pa_t w|^q) \|_{Z^2(\kappa_2)} \le
 C \|w\|_{Z^2(\kappa_1)}^q.
\end{align}
\end{lemma}

\noindent{\it Proof.}\ 
For $v \in X^2(\kappa_2)$, we have
$$
[v(r,t)]_2 \le \|v\|_{X^2(\kappa_2)}(1+|r-t|)^{-\kappa_2},
$$ 
so that $M_1(|\pa_t v|^p) \le p \|v\|_{X^2(\kappa_2)}^p$
holds for $\alpha=0$, $\beta=p$, $\gamma=0$, and
$\delta=p\kappa_2$.
By Proposition 3.1 with $s=2$ we get $L(|\pa_t v|^p) \in Z^2(\kappa_1)$
and \eqref{ta1}, since $\alpha+\beta+\gamma-1=\kappa_1<1$
and $\delta=p\kappa_2>1$.

On the other hand, for $w \in Z^2(\kappa_1)$, we have
$$
[w(r,t)]_2 \le \|w\|_{Z^2(\kappa_1)}(1+r+t)^{1-\kappa_1}
(1+|r-t|)^{-1},
$$ 
so that $M_1(|\pa_t w|^q) \le q \|w\|_{Z^2(\kappa_1)}^q$
holds for $\alpha=0$, $\beta=\delta=q$, and
$\gamma=q\kappa_1-q=\kappa_2+1-q$.
By Proposition 3.2 with $s=2$ we get $R(|\pa_t w|^q) \in Z^2(\kappa_2)$
and \eqref{ta2}, since $\alpha+\beta+\gamma-1=\kappa_2>1$.
This completes the proof.
\hfill$\square$ 

Next we examine the qualitative property of 
$\{w_j\}_{j=0}^{\ell+1}$ and $\{v_j\}_{j=0}^{\ell+1}$
defined by \eqref{1.25+} and \eqref{1.26+}.
As a corollary of Lemma \ref{Lemma1}, we derive the following estimates.

\begin{corollary} \label{Lemma3.0}
Let $0 \le j \le \ell+1$, $0<\ve \le 1$ and
$\vec{f}_i \in Y_{\kappa_i}(\ve)$ with $i=1, 2$.
Then $w_{j} \in Z^2(\kappa_1)$, $v_{j} \in X^2(\kappa_2)$,
and we have
\begin{align} \label{kubota4.3} 
& \|w_{j}\|_{Z^2(\kappa_1)} \le C\ve,
\\ \label{kubota4.4} 
& \|v_{j}\|_{X^2(\kappa_2)} \le C\ve.
\end{align}
Besides, we have 
\begin{align} \label{kubota4.5} 
& \|w_{1}-w_{0}\|_{Z^2(\kappa_1)} \le C\ve^p,
\\ \label{kubota4.6} 
& \|v_{1}-v_{0}\|_{Z^2(\kappa_2)} \le C\ve^q.
\end{align}
\end{corollary}

\noindent{\it Proof.}\ 
Since $\vec{f}_2 \in Y_{\kappa_2}(\ve)$, by Proposition \ref{Th.3.1} we get
$v_{0} \in X^2(\kappa_2)$ and (\ref{kubota4.4}) for $j=0$.
Analogously we have $w_{0} \in X^2(\kappa_1)$ and
$\|w_{0}\|_{X^2(\kappa_1)} \le C\ve$. Since $0<\kappa_1<1$, we find
$w_{0} \in Z^2(\kappa_1)$ and (\ref{kubota4.3}) for $j=0$.

Next suppose that \eqref{kubota4.3} and \eqref{kubota4.4} hold
for some $j$ with $0 \le j \le \ell$.
Since $w_{j+1}-w_0=L(|\pa_t v_j|^p)$ by \eqref{1.25+}, we have
$w_{j+1}-w_0 \in Z^2(\kappa_1)$ and $\|w_{j+1}-w_0\|_{Z^2(\kappa_1)}
 \le C \ve^p$, using \eqref{ta1} and \eqref{kubota4.4}.
Hence we get $ \|w_{j+1}\|_{Z^2(\kappa_1)} \le C\ve$ for $0<\ve \le 1$
and \eqref{kubota4.5} by taking $j=0$ in the above.
Analogously we obtain $\|v_{j+1}-v_0\|_{Z^2(\kappa_2)}
 \le C \ve^q$ by \eqref{1.26+}, \eqref{ta2} and \eqref{kubota4.3}.
Therefore, $ \|v_{j+1}\|_{X^2(\kappa_2)} \le C\ve$ for $0<\ve \le 1$
and \eqref{kubota4.6} holds, because $\kappa_2>1$.
The proof is complete.
\hfill$\square$

The following estimates are crucial in the proof of Theorem \ref{Th.2} 
for $\ell \ge 1$.

\begin{lemma} \label{Lemma3.1}
Let $1 \le j \le \ell$, $0<\ve \le 1$ and
$\vec{f}_i \in Y_{\kappa_i}(\ve)$ with $i=1, 2$.
Then $w_{j+1}-w_{j} \in Z^2(\kappa_1 a_j)$, 
$v_{j+1}-v_{j} \in Z^2(a_{j+1})$,
and we have
\begin{align} \label{kubota4.7} 
& \|w_{j+1}-w_{j}\|_{Z^1(\kappa_1 a_j)} \le C\ve^{b_{j-1}+p-1},
\\ \label{kubota4.8} 
& \|v_{j+1}-v_{j}\|_{Z^1(a_{j+1})} \le C\ve^{b_j},
\end{align}
where we put $b_k=q+k(p+q-2)$ for a nonnegative integer $k$,
and 
\begin{align} \label{kubota4.9}
& \|w_{j+1}-w_{j}\|_{Z^2(\kappa_1 a_j)} \le C\ve^{B_{j-1}},
\\
\label{kubota4.10} 
& \|v_{j+1}-v_{j}\|_{Z^2(a_{j+1})} \le C\ve^{B_{j-1}+q-1},
\end{align}
where we put $B_k=1+(p-1)b_k$ for a nonnegative integer $k$.
\end{lemma}

\noindent{\it Proof.}\ 
Observe that \eqref{kubota4.6} implies \eqref{kubota4.8} for $j=0$,
since $b_0=q$ and $a_1=\kappa_2$.

First we show that if \eqref{kubota4.8} holds for some $j$ 
with $0 \le j \le \ell-1$, then \eqref{kubota4.7} with $j$ replaced 
by $j+1$ holds.
It follows from \eqref{1.25+} that
\begin{align} \label{kubota4.11}
 w_{j+2}-w_{j+1}=L(G(v_{j+1},v_j)),
 \quad G(v,v^*)=|\pa_t v|^p-|\pa_t v^*|^p
\end{align}
for $0 \le j \le \ell-1$.
Note that if $v$, $v^*\in {X^2(\kappa_2)}$ and
$v-v^*\in {Z^1(a_{j+1})}$, then we have 
\begin{align} \label{kubota4.12} 
& \quad  |G(v,v^*)(r,t)| 
\\ \nonumber
& \le p|\pa_t ( v-v^*)|
    (|\pa_t v|+|\pa_t v^*|)^{p-1}
\\ \nonumber
 & \le p \|v-v^*\|_{Z^1(a_{j+1})}
 (\|v\|_{X^2(\kappa_2)}+\|v^*\|_{X^2(\kappa_2)})^{p-1}
\\ \nonumber 
 & \quad \ \times r^{-1}(1+r)^{-(p-1)}(1+r+t)^{-(a_{j+1}-1)}(1+|r-t|)^{-1-\kappa_1\kappa_2}.
\end{align}
In addition, we have
\begin{align} \label{kubota4.13} 
 & (1+r+t)^{-(a_{j+1}-1)}(1+|r-t|)^{-1-\kappa_1\kappa_2}
\\ \nonumber
 & \le (1+r+t)^{-\kappa_1(a_{j+1}-1)}(1+|r-t|)^{-\kappa_1-\kappa_2},
\end{align}
since $0<\kappa_1<1$ and $a_{j+1} \ge \kappa_2$.

Applying \eqref{kubota4.12} to $G(v_{j+1},v_j)$ and using \eqref{kubota4.4},
\eqref{kubota4.8} and \eqref{kubota4.13}, we obtain
\begin{align*} 
 M_0(G(v_{j+1},v_j)) \le C\ve^{b_j+p-1}
\end{align*}
for $\alpha=1$, $\beta=\kappa_1$, $\gamma=\kappa_1a_{j+1}-\kappa_1\,(>0)$,
and $\delta=\kappa_1+\kappa_2$.
Since $\alpha+\beta+\gamma-1=\kappa_1 a_{j+1}\le \kappa_1 a_{\ell}$ 
and $\kappa_1 a_\ell<1$ from the assumption in Theorem \ref{Th.2},
if we apply \eqref{3.7+} with $s=1$ to the right hand side on \eqref{kubota4.11}, 
then the desired estimate holds.
In particular, we have \eqref{kubota4.7} with $j=1$, since
\eqref{kubota4.8} is valid for $j=0$. 

Next we show that if \eqref{kubota4.7} holds for some $j$ 
with $1 \le j \le \ell$, then \eqref{kubota4.8} is valid for the same $j$.
It follows from \eqref{1.26+} that
\begin{align} \label{kubota4.14}
 v_{j+1}-v_{j}=R(H(w_{j+1},w_j)),
 \quad H(w,w^*)=|\pa_t w|^q-|\pa_t w^*|^q
\end{align}
for $1 \le j \le \ell$.
Note that if $w$, $w^*\in {Z^2(\kappa_1)}$ and
$w-w^*\in {Z^1(\kappa_1a_{j})}$, then we have 
\begin{align} \label{kubota4.15} 
& \quad  |H(w,w^*)(r,t)| 
\\ \nonumber
& \le q|\pa_t ( w-w^*)|
    (|\pa_t w|+|\pa_t w^*|)^{q-1}
\\ \nonumber
 & \le q \|w-w^*\|_{Z^1(\kappa_1 a_{j})}
 (\|w\|_{Z^2(\kappa_1)}+\|w^*\|_{Z^2(\kappa_1)})^{q-1}
\\ \nonumber 
 & \quad \ \times r^{-1}(1+r)^{-(q-1)}
(1+r+t)^{-(q-1)(\kappa_1-1)-(\kappa_1 a_{j}-1)}(1+|r-t|)^{-q}.
\end{align}
Applying \eqref{kubota4.15} to $H(w_{j+1},w_j)$ and using \eqref{kubota4.3},
\eqref{kubota4.7}, we obtain
\begin{align*} 
 M_0(H(w_{j+1},w_j)) \le C\ve^{b_j}
\end{align*}
for $\alpha=1$, $\beta=q-1$, $\gamma=a_{j+1}+1-q$, and $\delta=q$.
Since $\alpha+\beta+\gamma-1=a_{j+1}>1$, 
if we apply \eqref{3.16} with $s=1$ to the right hand side on \eqref{kubota4.14}, 
then the desired estimate holds.
In conclusion, we have proven \eqref{kubota4.7} and \eqref{kubota4.8}
for $1 \le j \le \ell$. 

Next we show \eqref{kubota4.9} and \eqref{kubota4.10}.
Observe that if we put $B_{-1}=1$, then \eqref{kubota4.10} with $j=0$ follows from
\eqref{kubota4.6}. 

First we show that if \eqref{kubota4.10} holds for some $j$ 
with $0 \le j \le \ell-1$, then it, in combination with
\eqref{kubota4.8}, implies \eqref{kubota4.9} with $j$ replaced 
by $j+1$.
Note that if $v$, $v^*\in {X^2(\kappa_2)}$ and
$v-v^*\in {Z^2(a_{j+1})}$, then we have 
\begin{align} \label{kubota4.16} 
& \quad  \{(1+r)|G(v,v^*)(r,t)| +r|\pa_r G(v,v^*)(r,t)|\}
\\ \nonumber
& \quad \ \times r^{p-1}(1+r+t)^{\kappa_1(a_{j+1}-1)}(1+|r-t|)^{\kappa_1+\kappa_2}
\\ \nonumber
& \le 2p \{ \|v-v^*\|_{Z^1(a_{j+1})}^{p-1} \|v\|_{X^2(\kappa_2)}
\\ \nonumber
& \quad \quad +
\|v-v^*\|_{Z^2(a_{j+1})}
 (\|v\|_{X^2(\kappa_2)}+\|v^*\|_{X^2(\kappa_2)})^{p-1}\}.
\end{align}
In fact, similarly to \eqref{kubota4.12}, we have
\begin{align*} 
& \quad  |G(v,v^*)(r,t)| 
\le p \|v-v^*\|_{Z^2(a_{j+1})}
 (\|v\|_{X^2(\kappa_2)}+\|v^*\|_{X^2(\kappa_2)})^{p-1}
\\ \nonumber 
 & \quad \hspace{28mm}
 \times (1+r)^{-p}(1+r+t)^{-(a_{j+1}-1)}(1+|r-t|)^{-1-\kappa_1\kappa_2}.
\end{align*}
Since $1<p<2$, we obtain
\begin{align*} 
 & |\pa_r G(v,v^*)(r,t)|
\\ \nonumber
\le & 2p|\pa_t ( v-v^*)|^{p-1}
    |\pa_r \pa_t v|
 +p|\pa_r \pa_t ( v- v^*)|
    |\pa_t v^*|^{p-1}
\\ \nonumber
\le & 2p \|v-v^*\|_{Z^1(a_{j+1})}^{p-1} \|v\|_{X^2(\kappa_2)}
\\ \nonumber
 & \quad \times r^{-p} (1+r+t)^{-\kappa_1(a_{j+1}-1)}
    (1+|r-t|)^{-\kappa_1-\kappa_2}
\\ \nonumber
 & \ +p \|v-v^*\|_{Z^2(a_{j+1})} \|v^*\|_{X^2(\kappa_2)}^{p-1}   
\\ \nonumber
 & \quad \times r^{-1} (1+r)^{-(p-1)}(1+r+t)^{-(a_{j+1}-1)}(1+|r-t|)^{-1-\kappa_1\kappa_2}.
\end{align*}
By \eqref{kubota4.13} we get \eqref{kubota4.16}.

Applying \eqref{kubota4.16} to $G(v_{j+1},v_j)$ and using \eqref{kubota4.4},
\eqref{kubota4.8} and \eqref{kubota4.10}, we obtain
\begin{align*} 
 M_1(G(v_{j+1},v_j)) \le C\ve^{(p-1)b_j+1}+C\ve^{B_{j-1}+p+q-2}
\end{align*}
for $\alpha=p-1\,(<1)$, $\beta=1$, $\gamma=\kappa_1a_{j+1}-\kappa_1$,
and $\delta=\kappa_1+\kappa_2$.
It is not difficult to see that
\begin{align} \label{kubota4.17}
B_j \le B_{j-1}+p+q-2
\end{align}
for $0 \le j \le \ell$ (recall $B_{-1}=1$).
Therefore we have $M_1(G(v_{j+1},v_j)) \le C\ve^{B_j}$.
Since $\alpha+\beta+\gamma-1=\kappa_1 a_{j+1}\le \kappa_1 a_{\ell} <1$,
if we apply \eqref{3.7+} with $s=2$ to the right hand side on \eqref{kubota4.11}, 
then the desired estimate holds.
In particular, we have \eqref{kubota4.9} with $j=1$, since
\eqref{kubota4.10} is valid for $j=0$. 

Finally we show that if \eqref{kubota4.9} holds for some $j$ 
with $1 \le j \le \ell$, then \eqref{kubota4.10} is valid for the same $j$.
Note that if $w$, $w^*\in {Z^2(\kappa_1)}$ and
$w-w^*\in {Z^2(\kappa_1a_{j})}$, then we have 
\begin{align} \label{kubota4.18} 
& \quad  \{(1+r)|H(w,w^*)(r,t)| +r|\pa_r H(w,w^*)(r,t)|\}
\\ \nonumber
& \quad \ \times (1+r)^{q-1}(1+r+t)^{(q-1)(\kappa_1-1)+\kappa_1 a_{j}-1}
   (1+|r-t|)^q
\\ \nonumber
& \le q^2 \|w-w^*\|_{Z^2(\kappa_1 a_{j})}
   (\|w\|_{Z^2(\kappa_1)}+\|w^*\|_{Z^2(\kappa_1)})^{q-1},
\end{align}
since $q>2$.
Applying \eqref{kubota4.18} to $H(w_{j+1},w_j)$ and using \eqref{kubota4.3},
\eqref{kubota4.9}, we obtain
\begin{align*} 
 M_1(H(w_{j+1},w_j)) \le C\ve^{B_{j-1}+q-1}
\end{align*}
for $\alpha=0$, $\beta=\delta=q$ and $\gamma=a_{j+1}+1-q$.
Since $\alpha+\beta+\gamma-1=a_{j+1}>1$, 
if we apply \eqref{3.16} with $s=2$ to the right hand side on \eqref{kubota4.14}, 
then the desired estimate holds.
In conclusion, all the asserion of the lemma is proven.
\hfill$\square$

Our next step is to solve the following system\,:
\begin{align}\label{kubota4.19} 
 u_1=w_{\ell+1}+R(G(u_2,v_{\ell})),
\quad u_2=v_{\ell+1}+R(H(u_1,w_{\ell+1})),
\end{align}
where $G(v,v^*)$ and $H(w,w^*)$ are the notations from \eqref{kubota4.11}
and \eqref{kubota4.14}, respectively.
We define $T(u_1,u_2)=(T^{(1)}(u_2),T^{(2)}(u_1))$ by
\begin{align}\label{kubota4.20} 
T^{(1)}(u_2)=w_{\ell+1}+R(G(u_2,v_{\ell})),
\quad T^{(2)}(u_1)=v_{\ell+1}+R(H(u_1,w_{\ell+1})).
\end{align} 
For $\ve>0$ we introduce a metric space
\begin{align} \label{kubota4.21} 
 D_\ve=\{(u_1,u_2) \in X^2 \times X^2\,;\,d((u_1,u_2),(w_{\ell+1},v_{\ell+1}))
   \le \ve^{(p-1)b_\ell} \},
\end{align}
where $b_\ell=q+\ell(p+q-2)$ and 
we have set $d((u_1,u_2),(u_1^*,u_2^*))=d_1(u_1,u_1^*) 
  +d_2(u_2,u_2^*)$ with
\begin{align}\nonumber
& d_1(u_1,u_1^*)=\| u_1-u_1^*\|_{Z^2(\kappa_1a_{\ell+1})}+
    \| u_1-u_1^*\|_{Z^1(\kappa_1a_{\ell+1})}^{p-1},
\\ \nonumber
& d_2(u_2,u_2^*)=\| u_2-u_2^*\|_{Z^2(a_{\ell+1})}+
   \| u_2-u_2^*\|_{Z^1(a_{\ell+1})}^{p-1}.
\end{align}
We shall show that $T$ is a contraction on $D_\ve$, provided
$\ve$ is small enough. 

First of all, we prepare the following.

\begin{lemma} \label{Lemma4.2}
Let $(u_1,u_2) \in D_\ve$ with $0<\ve \le 1$.
Then we have 
\begin{align} \label{kubota4.22} 
\|u_1\|_{Z^2(\kappa_1)}+\|u_2\|_{X^2(\kappa_2)} \le C\ve
\end{align}
and
\begin{align} \label{kubota4.23} 
d_1(u_1,w_{\ell+1}) \le \ve^{(p-1)b_\ell},
\quad d_2(u_2,v_{\ell}) \le C \ve^{(p-1)b_\ell}.
\end{align}
\end{lemma}

\noindent{\it Proof.}\ 
First we prove \eqref{kubota4.22}.
Notice that $a_{\ell+1}\ge \kappa_2>1$ and $(p-1)b_\ell \ge (p-1)q
=\kappa_2+1$.
Then \eqref{kubota4.22} follows from \eqref{kubota4.3} and 
\eqref{kubota4.4} with $j=\ell+1$.

Next we prove \eqref{kubota4.23}.
The first inequality is apparent.
On the other hand, in order to get the second one, it suffices
to show $d_2(v_{\ell+1},v_{\ell}) \le C \ve^{(p-1)b_\ell}$.
When $\ell=0$, it follows from \eqref{kubota4.6} that
$$
d_2(v_{1},v_{0}) \le C(\ve^q+\ve^{(p-1)q}) \le C\ve^{(p-1)q}
$$
for $0<\ve \le 1$.
While, when $\ell\ge 1$, it follows from \eqref{kubota4.8}
and \eqref{kubota4.10} with $j=\ell$ that
$$
d_2(v_{\ell+1},v_{\ell}) \le C(\ve^{B_{\ell-1}+q-1}+\ve^{(p-1)b_\ell}).
$$
Since 
\begin{align*} 
B_{\ell-1}+q-1=(p-1)b_\ell+(2-p)(p+q-1) >(p-1)b_\ell
\end{align*}
for $\ell \ge 1$, we obatin the needed estimate.
This completes the proof.
\hfill$\square$

The following estimate will play a basic role in proving
that $T$ is a contaction on $D_\ve$.

\begin{lemma} \label{Lemma4.5}
Let $u_2$, $u_2^* \in X^2(\kappa_2)$ satisfy
$u_2-u_2^* \in Z^2(a_{\ell+1})$ and
\begin{align} \label{kubota4.25} 
\|u_2\|_{X^2(\kappa_2)}+\|u_2^*\|_{X^2(\kappa_2)} \le C\ve.
\end{align}
Then we have
\begin{align}\label{kubota4.26} 
& \|R(G(u_2,u_2^*))\|_{Z^1(\kappa_1 a_{\ell+1})} \le 
 C \ve^{p-1} \|u_2-u_2^*\|_{Z^1(a_{\ell+1})} 
\\ \label{kubota4.27}
& \|R(G(u_2,u_2^*))\|_{Z^2(\kappa_1 a_{\ell+1})} \le 
 C \ve^{p-1} d_2(u_2,u_2^*).
\end{align}

While, let $u_1$, $u_1^* \in Z^2(\kappa_1)$ satisfy
$u_1-u_1^* \in Z^2(\kappa_1 a_{\ell+1})$ and
\begin{align} \label{kubota4.25+} 
\|u_1\|_{Z^2(\kappa_1)}+\|u_1^*\|_{Z^2(\kappa_1)} \le C\ve.
\end{align}
Then we have
\begin{align}\label{kubota4.33} 
& \|R(H(u_1,u_1^*))\|_{Z^1(a_{\ell+2})} \le 
 C \ve^{q-1} \|u_1-u_1^*\|_{Z^1(\kappa_1 a_{\ell+1})},
\\ \label{kubota4.34}
& \|R(H(u_1,u_1^*))\|_{Z^2(a_{\ell+2})} \le 
 C \ve^{q-1} \|u_1-u_1^*\|_{Z^2(\kappa_1 a_{\ell+1})}.
\end{align}
\end{lemma}

\noindent{\it Proof.}\ 
First we prove \eqref{kubota4.26}.
It follows from \eqref{kubota4.12} with $j=\ell$, 
\eqref{kubota4.13} and \eqref{kubota4.25} that
\begin{align*}
 M_0(G(u_{2},u_2^*)) \le C\ve^{p-1} \|u_2-u_2^*\|_{Z^1(a_{\ell+1})} 
\end{align*}
for $\alpha=1$, $\beta=\kappa_1$, $\gamma=\kappa_1a_{\ell+1}-\kappa_1$,
and $\delta=\kappa_1+\kappa_2$.
Applying \eqref{3.16} with $s=1$ to $G(u_2,u_2^*)$, we get \eqref{kubota4.26},
because 
\begin{align}\label{kubota4.28}
\alpha+\beta+\gamma-1=\kappa_1 a_{\ell+1}>1,
\end{align}
by virtue of \eqref{7}.

Next we prove \eqref{kubota4.27}.
It follows from \eqref{kubota4.16} with $j=\ell$, 
\eqref{kubota4.25} that
\begin{align*}
 M_1(G(u_{2},u_2^*)) \le C\ve^{p-1} d_2(u_2,u_2^*)
\end{align*}
for $\alpha=p-1$, $\beta=1$, $\gamma=\kappa_1a_{\ell+1}-\kappa_1$,
and $\delta=\kappa_1+\kappa_2$.
Since \eqref{kubota4.28} holds for these $\alpha$, $\beta$ and $\gamma$,
we obtain \eqref{kubota4.27} by \eqref{3.16} with $s=2$.

Next we prove \eqref{kubota4.33}.
It follows from \eqref{kubota4.15} with $j=\ell+1$
and \eqref{kubota4.25+} that
\begin{align*}
 M_0(H(u_{1},u_1^*)) \le C\ve^{q-1} 
 \|u_1-u_1^*\|_{Z^1(\kappa_1 a_{\ell+1})} 
\end{align*}
for $\alpha=1$, $\beta=q-1$, $\gamma=a_{\ell+2}+1-q$,
and $\delta=q$.
Since $\alpha+\beta+\gamma-1=a_{\ell+2}>1$, by \eqref{3.16} with $s=1$,
we get \eqref{kubota4.33}.

Finally we prove \eqref{kubota4.34}.
It follows from \eqref{kubota4.18} with $j=\ell+1$
and \eqref{kubota4.25+} that
\begin{align*}
 M_1(H(u_{1},u_1^*)) \le C\ve^{q-1} 
 \|u_1-u_1^*\|_{Z^2(\kappa_1 a_{\ell+1})} 
\end{align*}
for $\alpha=0$, $\beta=\delta=q$ and $\gamma=a_{\ell+2}+1-q$.
Since $\alpha+\beta+\gamma-1=a_{\ell+2}>1$, \eqref{3.16} with $s=2$
yields \eqref{kubota4.34}.
The proof is complete.
\hfill$\square$

\begin{corollary} \label{Cor4.7}
Let $(u_1,u_2) \in D_\ve$ with $0<\ve \le 1$.
Then we have 
\begin{align} \label{kubota4.29} 
& d_1(T^{(1)}(u_2), w_{\ell+1}) \le C\ve^{(p-1)b_\ell+(p-1)^2},
\\ \label{kubota4.31} 
& d_2(T^{(2)}(u_1),v_{\ell+1}) \le C\ve^{(p-1)b_\ell+(p-1)(q-1)}.
\end{align}
Moreover, we have
\begin{align} \label{kubota4.30} 
& d_1(T^{(1)}(u_2), T^{(1)}(u_2^*)) \le C\ve^{(p-1)^2}
   d_2(u_2,u_2^*),
\\ \label{kubota4.32} 
& d_2(T^{(2)}(u_1),T^{(2)}(u_1^*)) \le C\ve^{(p-1)(q-1)}
   d_1(u_1,u_1^*)
\end{align}
for $(u_1,u_2)$, $(u_1^*,u_2^*) \in D_\ve$ with $0<\ve \le 1$.
\end{corollary}

\noindent{\it Proof.}\ 
It follows from Lemma \ref{Lemma4.2} that
if $(u_1,u_2) \in D_\ve$ and $0<\ve \le 1$, then we have
\eqref{kubota4.22} and 
\begin{align}\label{kubota4.22+}
& \| u_1-w_{\ell+1}\|_{Z^2(\kappa_1a_{\ell+1})}+\| u_2-v_{\ell}\|_{Z^2(a_{\ell+1})}
   \le C \ve^{(p-1)b_\ell},
\\ \label{kubota4.22++}
& \| u_1-w_{\ell+1}\|_{Z^1(\kappa_1a_{\ell+1})}+\| u_2-v_{\ell}\|_{Z^1(a_{\ell+1})}
   \le C \ve^{b_\ell}.
\end{align}

We start with the proof of \eqref{kubota4.29}.
By \eqref{kubota4.20} we have $T^{(1)}(u_2)-w_{\ell+1}=R(G(u_2,v_{\ell}))$.
Therefore, applying the preceding lemma, we get \eqref{kubota4.29}.
Similarly, since $T^{(2)}(u_1)-v_{\ell+1}=R(H(u_1,w_{\ell+1}))$, 
we obtain \eqref{kubota4.31}.

Next we prove \eqref{kubota4.30}. By \eqref{kubota4.20} we have
$T^{(1)}(u_2)-T^{(1)}(u_2^*)=R(G(u_2,u_2^*))$.
Since $d_2(u_2,u_2^*) \le 2 \ve^{(p-1)b_\ell}$ for
$(u_1,u_2)$, $(u_1^*,u_2^*) \in D_\ve$, the preceding lemma
shows \eqref{kubota4.30}.
Similarly, we obtain \eqref{kubota4.32}. 
This completes the proof.
\hfill$\square$

\noindent
{\it End of the proof of Theorem \ref{Th.2}.}\
We see from Corollary \ref{Cor4.7} that
there exists a positive number $\ve_0$ depending only on
$p$ and $q$ such that if $0<\ve \le \ve_0$,
then we have $T(u_1,u_2) \in D_\ve$ and 
$$
d(T(u_1,u_2),T(u_1^*,u_2^*)) \le 2^{-1} 
d((u_1,u_2),(u_1^*,u_2^*))
$$ 
for $(u_1,u_2)$, $(u_1^*,u_2^*) \in D_\ve$, 
namely, $T$ is a contaction on $D_\ve$.
Hence we find a unique solution  
$(u_1,u_2) \in D_\ve$ of \eqref{kubota4.19}.
Here and in what follows, we suppose that $0<\ve \le \ve_0$
and $(u_1,u_2)$ is the solution.

Next we prove \eqref{12}.
Since \eqref{kubota4.4}, \eqref{kubota4.22} and \eqref{kubota4.22++} yield
\begin{align} \label{aux1}
\|u_2\|_{X^2(\kappa_2)}+\|v_\ell\|_{X^2(\kappa_2)} \le C\ve,
\quad
\|u_2-v_\ell\|_{Z^1(a_{\ell+1})} \le C\ve^{b_\ell},
\end{align}
applying \eqref{kubota4.26}, we obtain
\begin{align} \label{aux0} 
[R(G(u_2,v_\ell))(r,t)]_1 \le 
 C \ve^{p-1+b_\ell} (1+r+t)^{-(\kappa_1 a_{\ell+1}-1)}
 (1+|r-t|)^{-1}. 
\end{align} 
In view of \eqref{kubota4.19}, we have
\begin{align} \label{aux2}
\|(u_1-w_{\ell+1})(t)\|_E \le 
 C \ve^{p-1+b_\ell} (1+t)^{-(\kappa_1 a_{\ell+1}-1)}
\end{align} 
for $t \ge 0$.
Similarly, by \eqref{kubota4.3}, \eqref{kubota4.22}, \eqref{kubota4.22++}, and \eqref{kubota4.33},
 we get
\begin{align} \label{aux3}
\|(u_2-v_{\ell+1})(t)\|_E \le 
 C \ve^{q-1+b_\ell} (1+t)^{-(a_{\ell+2}-1)}
\end{align} 
for $t \ge 0$.
Moreover, it follows from \eqref{kubota4.6} and \eqref{kubota4.8} that
\begin{align} 
[(v_{\ell+1}-v_0)(r,t)]_1 \le 
 C \ve^{q} (1+r+t)^{-(\kappa_2-1)} (1+|r-t|)^{-1},
\end{align}
so that
\begin{align} \label{aux4} 
\|(v_{\ell+1}-v_0)(t)\|_E \le 
 C \ve^{q} (1+t)^{-(\kappa_{2}-1)}.
\end{align} 
Thus we obtain \eqref{12} from \eqref{aux2}, \eqref{aux3} and \eqref{aux4}.

Next we prove \eqref{13}.
Since $(w_{\ell+1}, \pa_t w_{\ell+1})(r,0)=\vec{f}_1(r)$,
it suffices to prove 
\begin{align} \label{kubota4.37} 
[R(G(u_2,v_\ell))(r,t)]_2 \le 
 C \ve^{B_\ell} (1+r+t)^{-(\kappa_1 a_{\ell+1}-1)}
 (1+|r-t|)^{-1}. 
\end{align} 
Using \eqref{kubota4.16} and \eqref{aux1}, we have
$$
M_1(G(u_2,v_\ell)) \le C(\ve^{B_\ell}+
   \ve^{p-1} \|u_2-v_\ell\|_{Z^2(a_{\ell+1})})
$$
for $\alpha=\kappa_1$, $\beta=1$, $\gamma=\kappa_1a_{\ell+1}-\kappa_1$,
and $\delta=\kappa_1+\kappa_2$.
It follows from \eqref{kubota4.10} and \eqref{kubota4.17} that 
$$
\ve^{p-1} \|v_{\ell+1}-v_\ell\|_{Z^2(a_{\ell+1})} \le C \ve^{B_\ell}
$$
for $0 <\ve \le 1$.
We see from \eqref{kubota4.19}, \eqref{kubota4.20} and \eqref{kubota4.31}
that
$$
\ve^{p-1} \|u_2-v_{\ell+1}\|_{Z^2(a_{\ell+1})} \le C \ve^{(p-1)b_\ell
  +(p-1)q}.
$$
Since $(p-1)q=\kappa_2+1>2$, we therefore obtain 
$$
M_1(G(u_2,v_\ell)) \le C \ve^{B_\ell}
$$
for $0 <\ve \le 1$.
Since \eqref{kubota4.28} is satisfied for those $\alpha$, $\beta$,
$\gamma$, and $\delta$, applying \eqref{3.16} with $s=2$ 
we get \eqref{kubota4.37}.

Finally we show \eqref{11}, which follows from
\begin{align} \label{aux5} 
\|u_2-v_0\|_{Z^2(\kappa_2)} \le  C \ve^q.
\end{align} 
Let $0<\ve \le 1$.
By \eqref{kubota4.18}, \eqref{kubota4.3} and \eqref{kubota4.22+}
we have
$$
M_1(H(u_1,w_{\ell+1})) \le C \ve^{q-1+(p-1)b_\ell} \le  C \ve^q
$$
for $\alpha=0$, $\beta=\delta=q$ and $\gamma=a_{\ell+2}+1-q$,
since $(p-1)b_\ell \ge (p-1)q >2$.
Noting $\alpha+\beta+\gamma-1=a_{\ell+2}>1$ and
applying \eqref{3.16} with $s=2$, we get
\begin{align*} 
 \|u_2-v_{\ell+1}\|_{Z^2(a_{\ell+2})} \le  C \ve^q,
\end{align*}
in view of \eqref{kubota4.19}.
While, we see from \eqref{kubota4.6} and \eqref{kubota4.10} that
\begin{align*} 
 \|v_{\ell+1}-v_0 \|_{Z^2(\kappa_2)} \le  C \ve^q,
\end{align*}
since $a_{j+1} \ge \kappa_2$ and $B_{j-1} \ge 1$ for $j \ge 1$.
Therefore we obtain \eqref{aux5}.
Thus we complete the proof of Theorem \ref{Th.2}.
\hfill $\square$

\vspace{2mm}

\noindent
{\bf Remark.}\ 
When $p=2$, Lemma \ref{Lemma1} and Corollary \ref{Lemma3.0} with $\ell=0$
remain valid, in view of \eqref{aux6}.
We replace \eqref{kubota4.21} by
\begin{align*} 
 D_\ve=\{(u_1,u_2) \in X^2 \times X^2\,;\,\|u_1-w_1\|_{Z^2(\kappa_2)}
+\|u_2-v_1\|_{Z^2(\kappa_2)} \le \ve^q\}.
\end{align*}
Notice that we have
$$ 
\|R(G(u_2,u_2^*))\|_{Z^2(\kappa_2)} \le C\ve \|u_2-u_2^*\|_{Z^2(\kappa_2)},
$$
instead of \eqref{kubota4.27}, because $M_1(G(u_2,u_2^*)) \le C\ve \|u_2-u_2^*\|_{Z^2(\kappa_2)}$
for $\alpha=0$, $\beta=2$, $\gamma=\kappa_2-1$, and $\delta=1+\kappa_2$
(remark that $\alpha+\beta+\gamma-1=\kappa_2>1$).
Prodceeding as in the proof of Theorem \ref{Th.2}, we find the desired conclusion stated 
after Theorem \ref{Th.2-}.

\vspace{2mm}

Next we prove Theorem \ref{Th.4}.
Similarly to the proofs of Corollary \ref{Lemma3.0} and Lemma \ref{Lemma3.1},
one can establish the following lemma.

\begin{lemma} \label{Lemma5.1}
Let $0<\ve \le 1$ and
$\vec{\varphi}_i \in Y_{\kappa_i}(\ve)$ with $i=1, 2$.
Then $w_{j}^* \in Z^2(\kappa_1)$, $v_{j}^* \in X^2(\kappa_2)$,
and we have
\begin{align} \label{kubota5.1} 
 \|w_{j}^* \|_{Z^2(\kappa_1)}+\|v_{j}^* \|_{X^2(\kappa_2)} \le C\ve
\end{align}
for $0 \le j \le \ell$. 
Moreover, $u_1-w_{j}^* \in Z^2(\kappa_1 a_j)$, 
$u_2-v_{j}^* \in Z^2(a_{j+1})$,
and we have
\begin{align} \label{kubota5.4} 
 \| u_1-w_{j}^* \|_{Z^1(\kappa_1 a_j)} \le C\ve^{b_{j-1}+p-1}, \quad
 \| u_2-v_{j}^*\|_{Z^1(a_{j+1})} \le C\ve^{b_j},
\end{align}
\begin{align} \label{kubota5.6}
 \| u_1-w_{j}^* \|_{Z^2(\kappa_1 a_j)} \le C\ve^{B_{j-1}}, \quad
 \| u_2-v_{j}^*\|_{Z^2(a_{j+1})} \le C\ve^{B_{j-1}+q-1}
\end{align}
for $1 \le j \le \ell$, together with
\begin{align} \label{kubota5.7} 
 \| u_2-v_{0}^*\|_{Z^2(\kappa_2)} \le C\ve^{q}.
\end{align}
Here $b_k$ and $B_k$ are defined in Lemma \ref{Lemma3.1}.
\end{lemma}

\noindent
{\it End of the proof of Theorem \ref{Th.4}.}\
Let $(u_1,u_2) \in Z^2(\kappa_1) \times {X^2(\kappa_2)}$ be the unique solution of 
\eqref{2.31} satisfying \eqref{1.36}.

In order to prove \eqref{2.40}, it suffices to show
\begin{align*} 
[R(G(u_2,v_\ell^*))(r,t)]_1 \le 
 C \ve^{p-1+b_\ell} (1+r+t)^{-(\kappa_1 a_{\ell+1}-1)}
 (1+|r-t|)^{-1},
\end{align*} 
in view of \eqref{2.34} and \eqref{2.38}.
The needed estimate can be deduced from \eqref{1.36}, \eqref{kubota5.1},
\eqref{kubota5.4}, and \eqref{kubota5.7}, similarly to \eqref{aux0}.

Next we show \eqref{1.39}.
By \eqref{2.34} and \eqref{2.35}, it is enough to prove
\begin{align*} 
[R(G(u_2,v_\ell^*))(r,t)]_2 \le 
 C \ve^{B_\ell} (1+r+t)^{-(\kappa_1 a_{\ell+1}-1)}
 (1+|r-t|)^{-1}.
\end{align*} 
Similarly to \eqref{kubota4.37}, we obtain the desired estimate
from \eqref{1.36}, \eqref{kubota5.1}, \eqref{kubota5.6}, and \eqref{kubota5.7}.

Finally we show \eqref{2.39} and \eqref{1.35}.
Since it follows from \eqref{2.38} that
$u_2-K[\vec{f}_2]=u_2-v_0^*$, we see that
\eqref{2.39} and \eqref{1.35} are consequences 
of \eqref{kubota5.7}.
This completes the proof of Theorem \ref{Th.4}.
\hfill $\square$

\vspace{2mm}

\section*{Acknowledgement} 
The research of the first author was partially supported by
Grant-in-Aid for Science Research (20224013), JSPS.

\end{document}